\documentclass[reqno]{amsart}
\usepackage{amsfonts}%
\usepackage{amssymb}
\usepackage{mathrsfs}%
\usepackage{cite}
\usepackage{color}
\allowdisplaybreaks
\date{\today}
\usepackage{mathrsfs}
\usepackage{amssymb}
\usepackage{amsmath}
\usepackage{amsthm}
\usepackage{latexsym}
\usepackage{graphicx}

 \numberwithin{equation}{section}

\newcommand{\dv}{\mathrm{div}\,}

\newcommand{\R}{{\mathbb{R}^3}}
\newcommand{\rt}{\mathrm{rot}\,}

\newtheorem{Theorem}{Theorem}[section]

\newtheorem{Remark}{Remark}[section]
\begin{document}

\title[Full Navier-Stokes-Maxwell system]
 {Convergence of the full compressible Navier-Stokes-Maxwell system
 to the incompressible magnetohydrodynamic equations in a bounded domain}

\author[J.-S. Fan]{Jishan Fan}
\address{ Department of Applied Mathematics,
 Nanjing Forestry University, Nanjing 210037, P.R.China}
\email{fanjishan@njfu.edu.cn}

\author[F.-C. Li]{Fucai Li$^*$}
\address{Department of Mathematics, Nanjing University, Nanjing
 210093, P.R. China}
 \email{fli@nju.edu.cn}
 \thanks{$^*$Corresponding author}
 \author[G. Nakamura]{Gen Nakamura}
\address{Department of Mathematics, Inha University, Incheon 402-751, Korea}
 \email{nakamuragenn@gmail.com}

%

\begin{abstract}
In this paper we  establish the uniform estimates of strong solutions with respect to
 the Mach number and the dielectric constant to the full compressible Navier-Stokes-Maxwell
system in a   bounded domain. Based on these uniform estimates, we obtain
the convergence of the  full compressible Navier-Stokes-Maxwell system
 to the incompressible magnetohydrodynamic equations  for well-prepared data.
 \end{abstract}


\maketitle

\section{Introduction}

In this paper we consider the singular limit of the following full compressible Navier-Stokes-Maxwell system in a bounded domain $
\Omega \subset \mathbb{R}^3$ (\!\cite{Im}):
\begin{align}
&\partial_t\rho+\dv(\rho u)=0,\label{1.1}\\
&\partial_t(\rho u)+\dv(\rho u\otimes u)+\frac{1}{\epsilon_1^2}\nabla p-\mu\Delta u-(\lambda+\mu)\nabla\dv u\nonumber\\
&\qquad\qquad=(E+u\times b)\times b,\label{1.2}\\
&\partial_t(\rho e)+\dv(\rho ue)+p\dv u-\dv(\kappa\nabla\mathcal{T})\nonumber\\
&\qquad\qquad=\epsilon_1^2(2\mu|D(u)|^2+\lambda(\dv u)^2+(E+u\times b)^2),\label{1.3}\\
&\epsilon_2\partial_tE-\rt b+E+u\times b=0,\label{1.4}\\
&\partial_tb+\rt E=0,\ \dv b=0,\label{1.5}
\end{align}
where the unknowns $\rho,u,p,e,\mathcal{T},E$, and $b$ stand for the density, velocity, pressure, internal energy, temperature, electric field, and magnetic field, respectively. The physical constants $\mu$ and $\lambda$ are the shear viscosity and bulk viscosity of the flow and satisfy $\mu>0$ and $\lambda+\frac23\mu\geq 0$. $\kappa>0$ is the heat conductivity. $\epsilon_1>0$ is the (scaled) Mach number, and $\epsilon_2>0$ is the (scaled) dielectric constant.
$D(u):=\frac12(\nabla u+\nabla u^\mathrm{T})$, where  $\nabla u^\mathrm{T}$ denotes the transpose
of the matrix $\nabla u$.

In \cite{KSa,KSb}, Kawashima and Shizuta   established  the global existence of smooth solutions for small data \cite{2}
 and studied its zero  dielectric constant  limit $\epsilon_2\rightarrow0$ in the whole space $\mathbb{R}^2$.
Recently, Jiang and Li \cite{JL} studied the  zero  dielectric constant  limit  $\epsilon_2\rightarrow0$ to the system \eqref{1.1}-\eqref{1.5}
 and obtained the convergence of the system \eqref{1.1}-\eqref{1.5} to the full compressible magnetohydrodynamic equations in $ \mathbb{T}^3$, see also
 \cite{JL2} on the similar results to the invisid case of \eqref{1.1}-\eqref{1.5}.
In \cite{LM}, Li and Mu study the low Mach number limit $\epsilon_1\rightarrow 0$  to the system \eqref{1.1}-\eqref{1.5}
and obtained the convergence of the system \eqref{1.1}-\eqref{1.5} to the incompressible Navier-Stokes-Maxwell system in the torus $\mathbb{T}^3$.

It should be pointed out that no boundary effect is  considered in the references mentioned above.
The purpose of this paper is to invistigate the singular limit $\epsilon_1,\epsilon_2\rightarrow 0$ to the system \eqref{1.5}-\eqref{1.5} in a bounded domain.
For simplicity, we shall take $\epsilon_1=\epsilon_2=\epsilon$ and consider the case that the fluid is a polytropic ideal gas, that is
\begin{equation}
e:=C_V\mathcal{T},\quad  p:=\mathcal{R}\rho\mathcal{T}\label{1.6}
\end{equation}
with $C_V>0$ and $\mathcal{R}$ being the specific heat at constant volume and the generic gas constant, respectively.

To state the main result of this paper, we denote the density and temperature variations by $\sigma^\epsilon$ and $\theta^\epsilon$:
\begin{equation}
\rho^\epsilon:=1+\epsilon\sigma^\epsilon,\ \mathcal{T}^\epsilon:=1+\epsilon\theta^\epsilon.\label{1.7}
\end{equation}
Then we can rewrite the  system \eqref{1.1}-\eqref{1.5} as follows:
\begin{align}
&\partial_t\sigma^\epsilon+\dv(\sigma^\epsilon u^\epsilon)+\frac{1}{\epsilon}\dv u^\epsilon=0,\label{1.8}\\
&\rho^\epsilon(\partial_tu^\epsilon+u^\epsilon\cdot\nabla u^\epsilon)+\frac{\mathcal{R}}{\epsilon}(\nabla\sigma^\epsilon+\nabla\theta^\epsilon)+\mathcal{R}\nabla(\sigma^\epsilon\theta^\epsilon)-\mu\Delta u^\epsilon-(\lambda+\mu)\nabla\dv u^\epsilon\nonumber\\
&\qquad\qquad=(E^\epsilon+u^\epsilon\times b^\epsilon)\times b^\epsilon,\label{1.9}\\
&C_V\rho^\epsilon(\partial_t\theta^\epsilon+u^\epsilon\cdot\nabla\theta^\epsilon)+\mathcal{R}(\rho^\epsilon\theta^\epsilon+\sigma^\epsilon)\dv u^\epsilon+\frac{\mathcal{R}}{\epsilon}\dv u^\epsilon\nonumber\\
&\qquad\qquad=\kappa\Delta\theta^\epsilon+\epsilon[2\mu|D(u^\epsilon)|^2+\lambda(\dv u^\epsilon)^2+(E^\epsilon+u^\epsilon\times b^\epsilon)^2],\label{1.10}\\
&\epsilon\partial_tE^\epsilon-\rt b^\epsilon+E^\epsilon+u^\epsilon\times b^\epsilon=0,\label{1.11}\\
&\partial_tb^\epsilon+\rt E^\epsilon=0,\ \dv b^\epsilon=0.\label{1.12}
\end{align}
Here we have added the superscript $\epsilon$ on the unknowns $(\sigma,u,\theta,E,b)$ to emphasise  the dependence of $\epsilon$.
The system \eqref{1.8}-\eqref{1.12} are supplemented with  the following initial and boundary conditions:
\begin{align}
&(\sigma^\epsilon,u^\epsilon,\theta^\epsilon,E^\epsilon,b^\epsilon)(\cdot,0)=(\sigma^\epsilon_0,u^\epsilon_0,\theta^\epsilon_0,E^\epsilon_0,b^\epsilon_0)(\cdot)\ \ \mathrm{in}\ \ \Omega,\label{1.13}\\
&u^\epsilon\cdot n=0,\ \rt u^\epsilon\times n=0,\ \frac{\partial\theta^\epsilon}{\partial n}=0,\ E^\epsilon\times n=0,\ b^\epsilon\cdot n=0\ \ \mathrm{on}\ \ \partial\Omega,\label{1.14}
\end{align}
where $n$ is the unit outer normal vector to the smooth boundary $\partial\Omega$.
%

Formally, if we let  $\epsilon\rightarrow 0$ in \eqref{1.8} and \eqref{1.9}, then we obtain that $\dv u^\epsilon \rightarrow 0, \nabla \theta^\epsilon \rightarrow 0,
$ and $\nabla \sigma^\epsilon \rightarrow 0$.  Letting  $\epsilon=0$ in \eqref{1.11} gives $E^\epsilon=\rt b^\epsilon-u^\epsilon\times b^\epsilon$. Pulling
it into   \eqref{1.12} and taking the limit  $\epsilon\rightarrow 0$ we obtain the following limit system (suppose that the limits $(u^\epsilon, b^\epsilon)\rightarrow (v, B)$ exist):
\begin{equation}
\left\{
\begin{array}{l}
v_t+v\cdot\nabla v+\nabla\pi-\mu\Delta v=\rt B\times B,\\
B_t+\rt(B\times v)-\Delta B=0,\\
\dv v=0,\quad \dv B=0.
\end{array}\right.\label{limit}
\end{equation}
We shall give a rigorous proof the the above formal analysis below.

Denote
\begin{align}
M^\epsilon(t):=&\sup\limits_{0\leq s\leq t}\Big\{\|(\sigma^\epsilon,u^\epsilon,\theta^\epsilon,\sqrt\epsilon E^\epsilon,b^\epsilon)(\cdot,s)\|_{H^2}
 +\|\partial_t(\sigma^\epsilon,u^\epsilon,\theta^\epsilon,\sqrt\epsilon E^\epsilon,b^\epsilon)(\cdot,s)\|_{H^1}\nonumber\\
&+\epsilon\|\partial_t^2(\sigma^\epsilon,u^\epsilon,\theta^\epsilon)(\cdot,s)\|_{L^2} +\Big\|\frac{1}{1+\epsilon\sigma^\epsilon(\cdot,s)}\Big\|_{L^\infty}\Big\}\nonumber\\
&+\bigg\{\int_0^t\Big(\|(u^\epsilon,\theta^\epsilon)\|_{H^3}^2+\|\partial_t(u^\epsilon,\theta^\epsilon)\|_{H^2}^2 +\|\epsilon\partial_t^2(\sigma^\epsilon,u^\epsilon,\theta^\epsilon)\|_{H^1}^2\nonumber\\
&+\|E^\epsilon\|_{H^2}^2+\|\partial_t(E^\epsilon,b^\epsilon)\|_{H^1}^2\Big)ds\bigg\}^\frac12.\label{1.16}
\end{align}

First, we have
\begin{Theorem}\label{th1.1}
Let $\Omega\subset \mathbb{R}^3$  be a simply connected, bounded domain with smooth boundary $\partial \Omega$ and $0<\epsilon<1$.
Suppose that the initial data $(\sigma_0^\epsilon,u_0^\epsilon,\theta_0^\epsilon,E_0^\epsilon,b_0^\epsilon)$ satisfy the following regularity conditions:
\begin{align}
&0\leq\theta_0^\epsilon,\quad 0<\frac{1}{K_0}\leq 1+\epsilon\sigma_0^\epsilon\leq K_0, \\
&\|(\sigma_0^\epsilon,u_0^\epsilon,\theta_0^\epsilon,E_0^\epsilon,b_0^\epsilon)\|_{H^2}+ \|\partial_t(\sigma^\epsilon,u^\epsilon,\theta^\epsilon,E^\epsilon,b^\epsilon)(\cdot,0)\|_{H^1}\nonumber\\
&\qquad\qquad \qquad \qquad \quad\ \, +\epsilon\|\partial_t^2(\sigma^\epsilon,u^\epsilon,\theta^\epsilon)(\cdot,0)\|_{L^2}\leq K_1\label{1.15}
\end{align}
for some positive constants $K_0>1$ and $K_1$ independent of $\epsilon>0$.
 Then there exist a small time $\tilde T>0$ independent of $\epsilon>0$ and a unique strong solution $(\sigma,u,\theta,E,b)$ to the initial boundary value problem \eqref{1.8}-\eqref{1.14} such that
\begin{equation}
M^\epsilon(\tilde T)\leq K\label{1.17}
\end{equation}
for some positive constant $K$ independent of $\epsilon>0$.
\end{Theorem}

\begin{Remark}
In the assumption \eqref{1.15}, $\sigma_t^\epsilon(\cdot,0)$ is indeed defined by $-\dv (\sigma_0^\epsilon u_0^\epsilon)+\frac{1}{\epsilon}\dv u_0^\epsilon$
 through the density equation and the other quantities  are defined by an analogous way.
\end{Remark}


Based on the uniform estimates of the solutions, we can prove the following convergence result by applying the Arzel\'{a}-Ascoli¡¯s theorem
in a standard way.

\begin{Theorem}\label{th1.3}
Let $(\sigma^{\epsilon}, u^{\epsilon}, \theta^{\epsilon},E^{\epsilon}, b^{\epsilon})$ be the solution of the problem \eqref{1.8}-\eqref{1.14} with
initial data $(\sigma^\epsilon_0,u^\epsilon_0,\theta^\epsilon_0,E^\epsilon_0,b^\epsilon_0)$ satisfying  the conditions in Theorem \ref{th1.1}.
Assume further that the initial data $(\sigma^\epsilon_0,u^\epsilon_0,\theta^\epsilon_0,E^\epsilon_0,b^\epsilon_0)$ satisfy that
\begin{align*}
&(\epsilon\sigma^{\epsilon}_0, u^{\epsilon}_0,\epsilon\theta^\epsilon_0,  b^{\epsilon}_0)\rightarrow(0, v_0,0, B_0)
\  \textrm{strongly in} \ H^s\   \text{for any} \  0\leq  s<2 \ \ \mathrm{as}\ \ \epsilon\rightarrow0,\\
&E^{\epsilon}_0\rightarrow \rt B_0-v_0\times B_0 \  \textrm{strongly in} \ H^s\   \text{for any} \  0\leq  s<1 \ \ \mathrm{as}\ \ \epsilon\rightarrow0.
\end{align*}
Then $(\epsilon\sigma^{\epsilon}, u^{\epsilon}, \epsilon\theta^{\epsilon}, b^{\epsilon})
\rightarrow(0, v,0, B)$ strongly in $L^\infty(0, \tilde T; H^1)$ and
$E^{\epsilon}\rightarrow B-v\times B$ strongly in $L^\infty(0, \tilde T; L^2)$ as $ \epsilon\rightarrow0$,
where $(v, B)$ satisfies \eqref{limit} with the following initial and boundary conditions:
\begin{equation}
\left\{
\begin{array}{l}
v\cdot n=B\cdot n=0,\ \rt v\times n=\rt B\times n=0\ \ on\ \ \partial\Omega\times(0, \tilde T],\\
(v, B)(\cdot, 0)=(v_0, B_0)(\cdot)\ \ in\ \ \Omega\subseteq\R.
\end{array}\right.\label{limit2}
\end{equation}
\end{Theorem}

The remainder of this paper is  devoted to the proof of Theorem \ref{th1.1} which will be given in next section.

\section{Proof of Theorem \ref{th1.1}}

In this section we shall prove Theorem \ref{th1.1} by combining the ideas developed in \cite{1,2,6,7}. First, by taking the
similar arguments to that \cite{1,2}, we know that in order to prove \eqref{1.17}, it suffices to show the following inequality
\begin{equation}
M^\epsilon(t)\leq C_0(M^\epsilon(0))\exp[t^\frac14C(M^\epsilon(t))]\label{1.18}
\end{equation}
for $\forall t\in[0,\tilde T]$ and some given positive nondecreasing continuous functions $C_0(\cdot)$ and $C(\cdot)$.

Below we shall   omit the spatial
domain $\Omega$ in the integrals and   drop the superscript ``$\epsilon$'' of $\rho^\epsilon,\sigma^\epsilon,u^\epsilon,\theta^\epsilon$, etc. for the sake of simplicity; moreover, we write $M^\epsilon(t)$ and $ M^\epsilon(0)$ as $M$ and $ M_0$, respectively.
Since the physical constants  $\kappa,C_V$,  and $\mathcal{R}$ do not bring any essential difficulties in
our arguments, we shall take $\kappa=C_V=\mathcal{R}=1$.

We will also use the following two inequalities:
\begin{align}
&\|u\|_{H^s(\Omega)}\leq C(\|\dv u\|_{H^{s-1}(\Omega)}+\|\rt u\|_{H^{s-1}(\Omega)}+\|u\|_{H^{s-1}(\Omega)}+\|u\cdot n\|_{H^{s-1/2}(\partial\Omega)}),\label{1.19}\\
&\|u\|_{H^s(\Omega)}\leq C(\|\dv u\|_{H^{s-1}(\Omega)}+\|\rt u\|_{H^{s-1}(\Omega)}+\|u\|_{H^{s-1}(\Omega)}+\|u\times n\|_{H^{s-1/2}(\partial\Omega)}),\quad\label{1.20}
\end{align}
for any $u\in H^s(\Omega) $ with  $s\geq1$, which were obtained in \cite{3} and \cite{4} respectively.


Because the local existence for the problem \eqref{1.8}-\eqref{1.14} with fixed $\epsilon>0$ is essential similar to that in  \cite{za}, we only need to prove \eqref{1.18}.
We will use the methods developed  in \cite{6,7}.

First, by the same calculations as that in \cite{6}, we get
\begin{align}
&\Big\|\frac{1}{\rho}(\cdot,t)\Big\|_{L^\infty}+\|\rho(\cdot,t)\|_{H^2}\leq C_0(M_0)\exp(C\sqrt tM),\label{2.1}\\
&\|\rho_t(\cdot,t)\|_{H^1}\leq C(M).\label{2.2}
\end{align}

Now we use the same method as that in \cite{7} to prove some a priori estimates on $(E,b)$.

Testing \eqref{1.11} and \eqref{1.12} by $E$ and $b$, respectively, and summing up the results, we see that
\begin{align*}
&\frac12\frac{d}{dt}\int(\epsilon E^2+b^2)dx+\int E^2dx=\int(b\times u)E dx\\
&\qquad\qquad\leq\|b\|_{L^2}\|E\|_{L^2}\|u\|_{L^\infty}\leq\frac12\int E^2dx+C(M).
\end{align*}
Integrating the above inequality over $(0,t)$, we find that
\begin{equation}
\int(\epsilon E^2+b^2)dx+\int_0^t\int E^2dx ds\leq C_0(M_0)+tC(M).\label{2.3}
\end{equation}

Using \eqref{1.14} and the formula
\begin{equation}
-(u\times b)\times n=(b\cdot n)u-(n\cdot u)b=0\ \ \mathrm{on}\ \ \partial\Omega,\label{2.4}
\end{equation}
we infer that
\begin{equation}
\rt b\times n=0\ \ \mathrm{on}\ \ \partial\Omega.\label{2.5}
\end{equation}

Taking $\rt$ to \eqref{1.11} and \eqref{1.12}, testing the results by $\rt E$ and $\rt b$, respectively, summing up the results, and using \eqref{2.5} and integration by parts, we have
\begin{align*}
&\frac12\frac{d}{dt}\int(\epsilon|\rt E|^2+|\rt b|^2)dx+\int|\rt E|^2dx\\
&\qquad=-\int\rt(u\times b)\cdot\rt E dx\\
&\qquad\leq\frac12\int|\rt E|^2dx+C\|u\|_{H^2}^2\|\rt b\|_{L^2}^2\\
&\qquad\leq\frac12\int|\rt E|^2dx+C(M).
\end{align*}

Integrating the above inequality over $(0,t)$, we have
\begin{equation}
\int(\epsilon|\rt E|^2+|\rt b|^2)dx+\int_0^t\int|\rt E|^2dx ds\leq C_0(M_0)+tC(M).\label{2.6}
\end{equation}

Taking $\dv$ to \eqref{1.11} and testing the result by $\dv E$, we infer that
\begin{align*}
&\frac{\epsilon}{2}\frac{d}{dt}\int(\dv E)^2dx+\int(\dv E)^2dx=\int\dv(b\times u)\dv E dx\\
&\qquad\leq\frac12\int(\dv E)^2dx+C(M).
\end{align*}

Integrating the above inequality over (0,t), we deduce that
\begin{equation}
\epsilon\int(\dv E)^2dx+\int_0^t\int(\dv E)^2dx ds\leq C_0(M_0)+tC(M).\label{2.7}
\end{equation}

Taking $\partial_t$ to \eqref{1.11} and \eqref{1.12}, testing the results by $E_t$ and $b_t$, respectively, summing up the results, we get
\begin{align*}
\frac12\frac{d}{dt}\int(\epsilon|E_t|^2+|b_t|^2)dx+\int|E_t|^2dx&=\int\partial_t(b\times u)\partial_tE dx\\
& \leq(\|b_t\|_{L^2}\|u\|_{L^\infty}+\|b\|_{L^\infty}\|u_t\|_{L^2})\|E_t\|_{L^2}\\
& \leq\frac12\int|E_t|^2dx+C(M).
\end{align*}

Integrating the above inequality over $(0,t)$, we get
\begin{equation}
\int(\epsilon|E_t|^2+|b_t|^2)dx+\int_0^t\int|E_t|^2dx ds\leq C_0(M_0)+tC(M).\label{2.8}
\end{equation}
\eqref{1.12} and \eqref{2.5} give the boundary condition
\begin{equation}
\rt^2E\times n=0\ \ \mathrm{on}\ \ \partial\Omega.\label{2.9}
\end{equation}

Taking $\rt^2$ to \eqref{1.11} and \eqref{1.12}, testing the results by $\rt^2E$ and $\rt^2b$, respectively, summing up the results, we derive that
\begin{align*}
&\frac12\frac{d}{dt}\int(\epsilon|\rt^2E|^2+|\rt^2b|^2)dx+\int|\rt^2E|^2dx\\
&\qquad=\int\rt^2(b\times u)\rt^2 E dx\\
&\qquad\leq\frac12\int|\rt^2E|^2dx+C\|b\|_{H^2}^2\|u\|_{H^2}^2\\
&\qquad\leq\frac12\int|\rt^2E|^2dx+C(M).
\end{align*}
Integrating the above inequality over $(0,t)$, we have
\begin{equation}
\int(\epsilon|\rt^2E|^2+|\rt^2b|^2)dx+\int_0^t\int|\rt^2E|^2dx ds\leq C_0(M_0)+tC(M).\label{2.10}
\end{equation}

Taking $\nabla\dv$ to \eqref{1.11}, testing the result by $\nabla\dv E$, we have
\begin{align*}
&\frac{\epsilon}{2}\frac{d}{dt}\int|\nabla\dv E|^2dx+\int|\nabla\dv E|^2dx=\int\nabla\dv(b\times u)\cdot\nabla\dv E dx\\
&\qquad\leq\frac12\int|\nabla\dv E|^2dx+C\|b\|_{H^2}^2\|u\|_{H^2}^2\\
&\qquad\leq\frac12\int|\nabla\dv E|^2dx+C(M).
\end{align*}
Integrating the above inequality over $(0,t)$, we obtain
\begin{equation}
\epsilon\int|\nabla\dv E|^2dx+\int_0^t\int|\nabla\dv E|^2dx ds\leq C_0(M_0)+tC(M).\label{2.11}
\end{equation}

Taking $\partial_t\rt$ to \eqref{1.11} and \eqref{1.12}, testing the results by $\partial_t\rt E$ and $\partial_t\rt b$, respectively, summing up the results, and using \eqref{2.5}, we obtain
\begin{align*}
&\frac12\frac{d}{dt}\int(\epsilon|\rt E_t|^2+|\rt b_t|^2)dx+\int|\rt E_t|^2dx\\
&\qquad=\int\rt(b_t\times u+b\times u_t)\rt E_t dx\\
&\qquad\leq C(\|b_t\|_{H^1}\|u\|_{H^2}+\|b\|_{H^2}\|u_t\|_{H^1})\|\rt E_t\|_{L^2}\\
&\qquad\leq\frac12\int|\rt E_t|^2dx+C(M).
\end{align*}
Integrating the above inequality over $(0,t)$, we obtain
\begin{equation}
\int(\epsilon|\rt E_t|^2+|\rt b_t|^2)dx+\int_0^t\int|\rt E_t|^2dx ds\leq C_0(M_0)+tC(M).\label{2.12}
\end{equation}

Applying $\partial_t\dv$ to \eqref{1.11}, testing the result by $\dv E_t$, we have
\begin{align*}
&\frac{\epsilon}{2}\frac{d}{dt}\int(\dv E_t)^2dx+\int(\dv E_t)^2dx\\
&\qquad=\int\dv(b_t\times u+b\times u_t)\dv E_t dx\\
&\qquad\leq C(\|b_t\|_{H^1}\|u\|_{H^2}+\|u_t\|_{H^1}\|b\|_{H^2})\|\dv E_t\|_{L^2}\\
&\qquad\leq\frac12\int(\dv E_t)^2dx+C(M).
\end{align*}
Integrating the above inequality over $(0,t)$, we have
\begin{equation}
\epsilon\int(\dv E_t)^2dx+\int_0^t\int(\dv E_t)^2dx ds\leq C_0(M_0)+tC(M).\label{2.13}
\end{equation}

Now we use the method in \cite{6} to prove some a priori estimates on $(\sigma,u,\theta)$.
Testing \eqref{1.8}, \eqref{1.9} and \eqref{1.10} by $\sigma,u$ and $\theta$, respectively, summing up the results, we obtain
\begin{align*}
&\frac12\frac{d}{dt}\int(\sigma^2+\rho u^2+\rho\theta^2)dx+\int(\mu|\nabla u|^2+(\lambda+\mu)(\dv u)^2+|\nabla\theta|^2)dx\\
=&\int\dv u\left(-\frac12\sigma^2-\rho\theta^2\right)dx+\epsilon\int\theta(2\mu|D(u)|^2+\lambda(\dv u)^2+(E+u\times b)^2)dx\\
\leq&\|\nabla u\|_{L^\infty}(\|\sigma\|_{L^2}^2+\|\rho\|_{L^\infty}\|\theta\|_{L^2}^2)\\
\leq &+C\|\nabla u\|_{L^\infty}(\|\nabla u\|_{L^2}+\|E\|_{L^4}^2+\|u\|_{L^\infty}^2\|b\|_{L^4}^2)\|\theta\|_{L^2}\\
\leq&\|\nabla u\|_{L^\infty}C(M)\leq\|u\|_{H^3}C(M).
\end{align*}
Integrating the above inequality over $(0,t)$, we obtain
\begin{equation}
\int(\sigma^2+\rho u^2+\rho\theta^2)dx+\int_0^t\int(|\nabla u|^2+|\nabla\theta|^2)dx ds\leq C_0(M_0)+\sqrt tC(M).\label{2.14}
\end{equation}

Applying $\partial_t$ to \eqref{1.8}, \eqref{1.9} and \eqref{1.10}, we see that
\begin{align}
&\partial_{tt}+\frac{1}{\epsilon}\dv u_t=-\dv(\sigma u)_t,\label{2.15}\\
&\rho(u_{tt}+u\cdot\nabla u_t)+\frac{1}{\epsilon}(\nabla\sigma_t+\nabla\theta_t)-\mu\Delta u_t-(\lambda+\mu)\nabla\dv u_t\nonumber\\
&\qquad=-\rho_tu_t-(\rho u)_t\nabla u-\nabla(\sigma\theta)_t+[(E+u\times b)\times b]_t,\label{2.16}\\
&\rho(\theta_{tt}+u\cdot\nabla\theta_t)+\frac{1}{\epsilon}\dv u_t-\Delta\theta_t=-\rho_t\theta_t-(\rho u)_t\cdot\nabla\theta-((\rho\theta+\sigma)\dv u)_t\nonumber\\
&\qquad+\epsilon(2\mu|D(u)|^2+\lambda(\dv u)^2+(E+u\times b)^2)_t.\label{2.17}
\end{align}
Testing \eqref{2.15}, \eqref{2.16} and \eqref{2.17} by $\sigma_t,u_t$ and $\theta_t$, respectively, summing up the results, we reach
\begin{align*}
&\frac12\frac{d}{dt}\int(\sigma_t^2+\rho u_t^2+\rho\theta_t^2)dx+\int(\mu|\nabla u_t|^2+(\lambda+\mu)(\dv u_t)^2+|\nabla\theta_t|^2)dx\\
&\qquad=\int(\sigma u)_t\nabla\sigma_t dx-\int[\rho_tu_t+(\rho u)_t\nabla u+\nabla(\sigma\theta)_t]u_t dx\\
&\qquad+\int((E+u\times b)\times b)_tu_t dx-\int[\rho_t\theta_t+(\rho u)_{t}\nabla\theta+((\rho\theta+\sigma)\dv u)_t]\theta_t dx\\
&\qquad+\epsilon\int(2\mu|D(u)|^2+\lambda(\dv u)^2+(E+u\times b)^2)_t\theta_t dx\\
&\qquad\leq C(M)+\int((E+u\times b)\times b)_tu_t dx+\epsilon\int((E+u\times b)^2)_t\theta_t dx\\
&\qquad\leq C(M)+\|E_t\|_{L^2}C(M).
\end{align*}
Integrating the above inequality over $(0,t)$, we have
\begin{equation}
\int(\sigma_t^2+\rho u_t^2+\rho\theta_t^2)dx+\int_0^t\int(|\nabla u_t|^2+|\nabla\theta_t|^2)dx ds\leq C_0(M_0)+\sqrt tC(M).\label{2.18}
\end{equation}

Testing \eqref{2.16} by $-\nabla\dv u$ in $L^2(\Omega\times(0,t))$, we find that
\begin{align}
&\frac{\mu+\lambda}{2}\|\nabla\dv u(\cdot,t)\|_{L^2}^2-\frac{1}{\epsilon}\int_0^t\int(\nabla\sigma_t+\nabla\theta_t)\cdot\nabla\dv u dx ds\nonumber\\
=&\frac{\mu+\lambda}{2}\|\nabla\dv u_0\|_{L^2}^2+\int_0^t\int\rho(u_{tt}+u\cdot\nabla u_t)\nabla\dv u dx ds\nonumber\\
&+\int_0^t\int(\rho_tu_t+(\rho u)_t\nabla u+\nabla(\sigma\theta)_t)\cdot\nabla\dv u dx ds\nonumber\\
&-\int_0^t\int((E+u\times b)\times b)_t\nabla\dv u dx ds\nonumber\\
=&:\frac{\mu+\lambda}{2}\|\nabla\dv u_0\|_{L^2}^2+I_1+I_2+I_3.\label{2.19}
\end{align}

We bound $I_1,I_2$ and $I_3$ as follows.
\begin{align*}
I_1\leq&C(M)\int_0^t\|u_{tt}\|_{L^2}ds+tC(M)\leq\sqrt tC(M),\\
I_2\leq&tC(M),\\
I_3\leq&C(M)\int_0^t\|E_t\|_{H^1}ds+tC(M)\leq\sqrt tC(M).
\end{align*}

Applying $\nabla$ to \eqref{1.8} and \eqref{1.10}, testing the results by $\nabla\sigma_t$ and $\nabla\theta_t$ in $L^2(\Omega\times(0,t))$, respectively, we derive
\begin{align}
&\int_0^t\int|\nabla\sigma_t|^2dx ds+\frac1\epsilon\int_0^t\int\nabla\sigma_t\nabla\dv u dx ds\nonumber\\
&\qquad=-\int_0^t\int\nabla\dv(\sigma u)\nabla\sigma_tdx ds\leq tC(M),\label{2.20}
\end{align}
and
\begin{align}
&\frac12\int|\Delta\theta|^2dx+\int_0^t\int\rho|\nabla\theta_t|^2dx ds+\frac1\epsilon\int_0^t\int\nabla\theta_t\nabla\dv u dx ds\nonumber\\
&\qquad=\frac12\int|\Delta\theta_0|^2dx+\int_0^t\int\epsilon\nabla[2\mu|D(u)|^2+\lambda(\dv u)^2+(E+u\times b)^2]\nabla\theta_tdx ds\nonumber\\
&\qquad\quad-\int_0^t\int\nabla[(\rho\theta+\sigma)\dv u]\nabla\theta_tdx ds-\int_0^t\int[\nabla\rho\theta_t+\nabla(\rho u\cdot\nabla\theta)]\nabla\theta_tdx ds\nonumber\\
&\qquad\leq C_0(M_0)+\sqrt tC(M).\label{2.21}
\end{align}

Summing up \eqref{2.19}, \eqref{2.20} and \eqref{2.21}, we arrive at
\begin{equation}
\int(|\nabla\dv u|^2+|\Delta\theta|^2)dx+\int_0^t\int(|\nabla\sigma_t|^2+|\nabla\theta_t|^2)dx ds\leq C_0(M_0)\exp(\sqrt tC(M)).\label{2.22}
\end{equation}

Testing \eqref{2.15}, \eqref{2.16} and \eqref{2.17} by $-\Delta\sigma_t,-\nabla\dv u_t$ and $-\Delta\theta_t$, respectively, we derive
\begin{align}
&\frac12\int|\nabla\sigma_t|^2dx+\frac12\int_0^t\int\nabla\dv u_t\nabla\sigma_tdx ds\nonumber\\
&\qquad=\frac12\int|\nabla\sigma_t(0)|^2dx+\int_0^t\int\dv(\sigma_tu+\sigma u_t)\cdot\Delta\sigma_tdx ds\nonumber\\
&\qquad=:\frac12\int|\nabla\sigma_t(0)|^2dx+I_4.\label{2.23}
\end{align}

We bound $I_4$ as follows.
\begin{align*}
I_4=&\int_0^t\int u\nabla\sigma_t\Delta\sigma_t dx-\int_0^t\int\nabla(\sigma_t\dv u+u_t\nabla\sigma+\sigma\dv u_t)\nabla\sigma_t dx ds\\
=&-\sum\limits_i\int_0^t\int\partial_iu\nabla\sigma_t\partial_i\sigma_tdx ds+\int_0^t\int\frac12\dv u|\nabla\sigma_t|^2dx ds\\
&-\int_0^t\int\nabla(\sigma_t\dv u+u_t\nabla\sigma+\sigma\dv u_t)\nabla\sigma_tdx ds\\
\leq&C(M)\int_0^t\|\nabla u\|_{L^\infty}ds+C(M)\int_0^t\|u\|_{H^3}ds+C(M)\int_0^t\|u_t\|_{H^2}ds\\
\leq&\sqrt tC(M).
\end{align*}

And
\begin{align}
&\frac12\int\rho(\dv u_t)^2dx+(\lambda+2\mu)\int_0^t\int|\nabla\dv u_t|^2dx ds\nonumber\\
&-\frac1\epsilon\int_0^t\int\nabla\dv u_t(\nabla\sigma_t+\nabla\theta_t)dx ds\nonumber\\
=&\frac12\int\rho_0(\dv u_t(0))^2dx+\int_0^t\int[\rho_tu_t+(\rho u)_t\nabla u+\nabla(\sigma\theta)_t]\nabla\dv u_tdxds\nonumber\\
&+\int_0^t\int\dv(E\times b)_t\dv u_tdxds-\int_0^t\int[(u\times b)\times b]_t\nabla\dv u_t dxds\nonumber\\
&+\int_0^t\int(\epsilon\nabla\sigma\cdot u_{tt}+u\cdot\nabla u_t)\dv u_t dxds-\int_0^t\int\sum\limits_i\nabla u_i\partial_iu_t\dv u_t dx ds\nonumber\\
\leq&C_0(M_0)+\sqrt tC(M).\label{2.24}
\end{align}

And
\begin{align}
&\frac12\int\rho|\nabla\theta_t|^2dx+\int_0^t\int|\Delta\theta_t|^2dx ds+\frac1\epsilon\int_0^t\int\nabla\dv u_t\nabla\theta_t dx ds\nonumber\\
=&\frac12\int\rho_0|\nabla\theta_t(0)|^2dx-\int_0^t\int\epsilon\nabla\sigma\theta_{tt}\nabla\theta_tdx ds-\int_0^t\int\sum\limits_i\nabla(\rho u_i)\partial_i\theta_t\nabla\theta_tdx ds\nonumber\\
&+\int_0^t\int\Delta\theta[\rho_t\theta_t+(\rho u)_t\nabla\theta+((\rho\theta+\sigma)\cdot\dv u)_t]dx ds\nonumber\\
&-\epsilon\int_0^t\int\Delta\theta_t(2\mu|D(u)|^2+\lambda(\dv u)^2+(E+u\times b)^2)_tdx ds\nonumber\\
\leq&C_0(M_0)+\sqrt tC(M).\label{2.25}
\end{align}

Summing up \eqref{2.23}, \eqref{2.24} and \eqref{2.25}, we arrive at
\begin{align}
&\int(|\nabla\sigma_t|^2+(\dv u_t)^2+|\nabla\theta_t|^2)dx+\int_0^t\int(|\nabla\dv u_t|^2+(\Delta\theta_t)^2)dx ds\nonumber\\
&\qquad\leq C_0(M_0)\exp(\sqrt tC(M)).\label{2.26}
\end{align}

Now, testing $\partial_i\nabla$ \eqref{1.8} by $\partial_i\nabla\sigma+\partial_i\nabla\theta$ and the same calculations as those in \cite{6} to obtain
\begin{align}
&\frac12\int|\partial_i\nabla\sigma|^2dx+\int\partial_i\nabla\sigma\cdot\partial_i\nabla\theta dx+\frac1\epsilon\int_0^t\int\partial_i\nabla\dv u(\partial_i\nabla\sigma+\partial_i\nabla\theta)dx ds\nonumber\\
&\qquad\leq C_0(M_0)+\sqrt tC(M).\label{2.27}
\end{align}

Testing $\partial_i$ \eqref{1.9} by $\partial_i\nabla\dv u$ in $L^2(\Omega\times(0,t))$ and the same calculations as those in \cite{6} to obtain
\begin{equation}
\frac12\int_0^t\int|\partial_i\nabla\dv u|^2dx-\frac1\epsilon\int_0^t\int\partial_i\nabla\dv u(\partial_i\nabla\sigma+\partial_i\nabla\theta)dx ds\leq tC(M).\label{2.28}
\end{equation}

\eqref{2.27}, \eqref{2.28} and \eqref{2.22} give
\begin{equation}
\int|\nabla^2\sigma|^2dx+\int_0^t\int|\nabla^2\dv u|^2dx ds\leq C_0(M_0)\exp(\sqrt tC(M)).\label{2.29}
\end{equation}

Applying $\rt$ to \eqref{1.9} and denoting the vorticity $\omega:=\rt u$, we see that
\begin{align}
\rho\omega_t+\rho u\cdot\nabla\omega-\mu\Delta\omega=&(\partial_j\rho u_{it}-\partial_i\rho u_{jt})+(\partial_j(\rho u_k)\partial_ku_i-\partial_i(\rho u_k)\partial_ku_j)\nonumber\\
& +\rt[(E+u\times b)\times b].\label{2.30}
\end{align}

We test \eqref{2.30} by $\Delta\omega$ in $L^2(\Omega\times(0,t))$ to get
\begin{equation}
\int\rho|\rt\omega|^2dx+\mu\int_0^t\int|\Delta\omega|^2dx ds\leq C_0(M_0)+\sqrt tC(M).\label{2.31}
\end{equation}

Similarly, we apply $\partial_t$ to \eqref{2.30} and test the resulting equations by $\omega_t$ in $L^2(\Omega\times(0,t))$ to deduce that
\begin{equation}
\int\rho|\omega_t|^2dx+\mu\int_0^t\int|\rt\omega_t|^2dx ds\leq C_0(M_0)+\sqrt tC(M).\label{2.32}
\end{equation}

By the same calculations as that in \cite{6}, we have
\begin{equation}
\|\Delta\theta\|_{L^2(0,t;H^1)}\leq C_0(M_0)\exp(t^\frac14C(M)).\label{2.33}
\end{equation}

It follows from \eqref{1.11}, \eqref{1.12}, \eqref{2.3}, \eqref{2.6}, \eqref{2.7}, \eqref{2.8}, \eqref{2.10}, \eqref{2.11}, \eqref{2.12} and \eqref{2.13} that
\begin{align}
&\epsilon\|E_{tt}\|_{L^2(0,t;L^2)}\leq\|\rt b_t-E_t-(u\times b)_t\|_{L^2(0,t;L^2)}\leq C_0(M_0)+tC(M).\label{2.34}\\
&\sqrt\epsilon\|b_{tt}(t)\|_{L^2}\leq\sqrt\epsilon\|\rt E_t(\cdot,t)\|_{L^2}\leq C_0(M_0)+tC(M).\label{2.35}
\end{align}

Finally, we need to estimate $\epsilon\sigma_{tt},\ \epsilon u_{tt}$ and $\epsilon\theta_{tt}$  to close the energy estimates. Testing $\partial_t^2$\eqref{1.8}, $\partial_t^2$\eqref{1.9} and $\partial_t^2$\eqref{1.10} by $\epsilon^2\sigma_{tt},\ \epsilon^2u_{tt}$ and $\epsilon^2\theta_{tt}$, respectively, then by the same calculations as that in \cite{6}, we conclude that $$\epsilon\|(\sigma_{tt},u_{tt},\theta_{tt})(t)\|_{L^2}+\epsilon\|(u_{tt},\theta_{tt})\|_{L^2(0,t;H^1)}\leq C_0(M_0)\exp(t^\frac14C(M))$$ and thus \eqref{1.17} hold true.

This completes the proof the proof of Theorem \ref{th1.1}.

\hfill$\square$




{\bf Acknowledgements:}
 Fan is supported by NSFC (Grant No. 11171154). Li is supported partially by NSFC (Grant No. 11271184) and
   PAPD.

\end{document}